 \documentclass[10pt]{article}
 \usepackage[top=1in,bottom=1in,left=1.5in,right=1.5in]{geometry}
 \usepackage{amsmath}
  \usepackage[dvips]{pstricks} 
  \usepackage{pst-node}
  \usepackage[dvips]{graphicx}
  \usepackage[standard,amsmath,thmmarks]{ntheorem}
  \usepackage{mathrsfs}

  \newcommand{\R}{\mathbb{R}}

  \newcommand{\cM}{\mathcal{M}}

  \newcommand{\cS}{\mathcal{S}}

  \newcommand{\lan}{\langle}
  \newcommand{\ran}{\rangle}
  \newcommand{\an}[1]{\lan#1\ran}
  \newcommand{\hs}{\hspace*{\parindent}}

  \newcommand{\trans}{^\top}
  \newcommand{\per}{\mathop{\mathrm{perm}}\nolimits}
  \newcommand{\permat}{\mathop{\mathrm{perfmat}}\nolimits}

  \newcommand{\haf}{\mathop{\mathrm{hafn}}\nolimits}

  \newcommand{\rS}{\mathrm{S}}

  \newtheorem{theo}{\bfseries \hs Theorem}[section]

  \newtheorem{prop}[theo]{\bfseries \hs Proposition}
  \newtheorem{lem}[theo]{\bfseries \hs Lemma}

  \numberwithin{equation}{section} 

\begin{document}
  \title{An upper bound for the number of\\ perfect matchings in graphs }
  \author{
  Shmuel Friedland\footnote{
   Visiting Professor, Fall 2007 - Winter 2008,
  Berlin Mathematical School, Berlin, Germany}
  }
  \date{Department of Mathematics, Statistics, and Computer Science,\\
        University of Illinois at Chicago\\
        Chicago, Illinois 60607-7045, USA\\
        6 March, 2008}

 \maketitle

 \begin{abstract}
 We give an upper bound on the number of perfect
 matchings in an undirected simple graph $G$ with an even number of
 vertices, in terms of the degrees of all the vertices in $G$.
 This bound is sharp if $G$ is a union of complete
 bipartite graphs.
 This bound is a generalization of the upper bound
 on the number of perfect matchings in bipartite graphs on
 $n+n$ vertices given by the Bregman-Minc inequality for the
 permanents of $(0,1)$ matrices.
     \\[\baselineskip] 2000 Mathematics Subject
     Classification: 05A15, 05C70.
 \end{abstract}


\section{Introduction}

 Let $G=(V,E)$ be an undirected simple graph with the set of vertices $V$
 and edges $E$. For a vertex $v\in V$ denote by $\deg v$ the
 degree of the vertex $v$.
 Assume that $\#V$ is even.  Denote by $\permat G$
 the number of perfect matching in $G$.  Our main result states
 that
 \begin{equation}\label{uppmatbd}
 \permat G\le \prod_{v\in V} ((\deg v)!)^{\frac{1}{2\deg v}},
 \end{equation}
 We assume here that $0^{\frac{1}{0}}=0$.
 This result is sharp if $G$ is a disjoint union of complete
 bipartite graphs.
 For bipartite graphs the above inequality follows from the
 Bregman-Minc inequality for the permanents of $(0,1)$ matrices, conjectured
 by Minc \cite{Min63} and proved by Bregman \cite{Bre}.
 In fact, the inequality (\ref{uppmatbd}) is the analog of the
 Bregman-Minc inequality for the \emph{hafnians} of $(0,1)$
 symmetric of even order with zero diagonal.
 Our proof follows closely the proof of the Bregman-Minc
 inequality given by Schrijver \cite{Sch78}.

 \section{Permanents and Hafnians}
 If $G$ is a bipartite
 graph on $n+n$ vertices then $\permat G= \per B(G)$, where $B(G)=[b_{ij}]\in
 \{0,1\}^{n\times n}$ is the incidence matrix of the bipartite graph $G$.
 Thus $V=V_1\cup V_2$ and $E\subset V_1\times V_2$, where
 $V_i=\{v_{1,i},\ldots,v_{n,i}\}$ for $i=1,2$.  Then $b_{ij}=1$
 if and only if $(v_{i,1},v_{j,2})\in E$.
 Recall that the permanent
 of $B\in \R^{n\times n}$ is given by
 $\per B=\sum_{\sigma\in\cS_n}
 \prod_{i=1}^n b_{i\sigma(i)}$,
 where $\cS_n$ is the symmetric group of
 all permutations $\sigma:\an{n}\to\an{n}$.

 Vice versa, given any $(0,1)$ matrix
 $B=[a_{ij}]\in\{0,1\}^{n\times n}$, then $B$ is the incidence
 matrix of the induced $G(B)=(V_1\cup V_2,E)$.
 Denote by $\an{n}:=\{1,\ldots,n\},
 m+\an{n}:=\{m+1,\ldots,m+n\}$ for any two positive integers
 $m,n$.  It is convenient to identify $V_1=\an{n},
 V_2=n+\an{n}$.
 Then $r_i:=\sum_{j=1}^n b_{ij}$ is the $i-th$
 degree of $i\in\an{n}$.
 The celebrated Bregman-Minc
 inequality, conjectured by Minc \cite{Min63} and proved by Bregman \cite{Bre},
 states
 \begin{equation}\label{bminq}
 \per B\le \prod_{i=1}^n (r_i!)^{\frac{1}{r_i}}.
 \end{equation}
 A simple proof Bregman-Minc inequality is given  \cite{Sch78}.
 Furthermore the above inequality is generalized to nonnegative
 matrices.  See \cite{AS, Rad} for additional proofs of
 (\ref{bminq}).

 \begin{prop}\label{upestbipgr}  Let $G=(V_1\cup V_2,E)$ be a
 bipartite graph with $\#V_1=\#V_2$.  Then (\ref{uppmatbd}) holds.
 If $G$ is a union of complete bipartite graphs then equality holds in
 (\ref{uppmatbd}).
 \end{prop}
 \begin{proof}
 Assume that $\#V_1=\#V_2=n$.
 Clearly,
 $$\permat G=\per B(G)=\per B(G)\trans=\sqrt{\per B(G)}\sqrt{\per
 B(G)\trans}.$$
 Note that the $i-th$ row sum of $B(G)\trans $ is the degree of
 the vertex $n+i\in V_2$.  Apply the Bregman-Minc inequality
 to $\per B(G)$ and $\per B(G)\trans$ to deduce
 (\ref{uppmatbd}).

 Assume that $G$ is the complete bipartite graph $K_{r,r}$ on
 $r+r$ vertices.  Then $B(K_{r,r})=J_r=\{1\}^{r\times r}$.
 So $\permat K_{r,r}=r!$.  Hence equality holds in
 (\ref{uppmatbd}).  Assume that $G$ is a (disjoint) union of $G_1,\ldots
 G_L$.
 Since $\permat G=\prod_{i=1}^L \permat G_i$, we deduce (\ref{uppmatbd})
 is sharp if each $G_i$ is a complete bipartite graph.
 \end{proof}

 Let $A(G)\in \{0,1\}^{m\times m}$ be the adjacency
 matrix of an undirected simple graph $G$ on $m$ vertices.
 Note that $A(G)$ is a symmetric matrix with zero diagonal.
 Vice versa, any symmetric $(0,1)$ matrix with zero diagonal
 induces an indirected simple graph $G(A)=(V,E)$ on $m$
 vertices.  Identify $V$ with $\an{m}$.
 Then $r_i$, the $i-th$ row sum of $A$, is the
 degree of the vertex $i\in\an{m}$.

  Let $K_{2n}$ be the complete graph on $2n$ vertices,
 and denote by $\cM(K_{2n})$ the set of all perfect matches
 in $K_{2n}$.  Then $\alpha\in \cM(K_{2n})$ can be represented as
 $\alpha=\{(i_1,j_1),(i_2,j_2),..,(i_n,j_n)\} $
 with $ i_k<j_k$ for $k\in\an{n}$.
 It is convenient to view $(i_k,j_k)$ as an edge in $K_{2n}$.
 We can view $\alpha$ as an involution in $\cS_{2n}$
 with no fixed points.  So for $l\in \an{2n}$ $\alpha(l)$ is
 second vertex corresponding to $l$ in the perfect match given
 by $\alpha$.  Vice versa, any fixed point free involution of
 $\an{2n}$ induces a perfect match $\alpha\in \cM(K_{2n})$.
 Denote by $\rS_m$ the space
 of $m\times m$ real symmetric matrices.  Assume that
 $A=[a_{ij}]\in \rS_{2n}$.  Then the \emph{hafnian} of $A$ is
 defined as
 \begin{equation}\label{defhaf}
 \haf A:=\sum_{\alpha=\{(i_1,j_1),(i_2,j_2),..,(i_n,j_n)\}\in\cM(K_{2n})}
 \prod_{k=1}^n a_{i_k j_k}.
 \end{equation}
 Note that $\haf A$ does not depend on the diagonal entries of $A$.
 Let $i\ne j\in\an{2n}$.  Denote by $A(i,j)\in \rS_{2n-2}$
 the symmetric matrix obtained from $A$ by deleting the $i,j$
 rows and columns of $A$.  The following proposition is
 straightforward, and is known as the expansion of the hafnian
 by the row, (column), $i$.

 \begin{prop}\label{hafrowex}  Let $A\in \rS_{2n}$.  Then for
 each $i\in\an{2n}$
 \begin{equation}\label{hafrowex1}
 \haf A= \sum_{j\in\an{2n}\backslash\{i\}} a_{ij} \haf A(i,j)
 \end{equation}
 \end{prop}

 It is clear that  $\permat G=\haf A(G)$ for any $G=(\an{2n},E)$.
 Then (\ref{uppmatbd}) is equivalent to the inequality
 \begin{equation}\label{bmhafin}
 \haf A\le \prod_{i=1}^{2n} (r_i!)^{\frac{1}{2r_i}} \textrm{ for all }
 A\in \{0,1\}^{(2n)\times (2n)}\cap \rS_{2n,0}
 \end{equation}

 Our proof of the above inequality follows the proof of the Bregman-Minc inequality
 given by A. Schrijver \cite{Sch78}.

 \section{Preliminaries}

 Recall that $x\log x$ is a strict convex function on
 $\R_+=[0,\infty)$, where $0\log 0=0$.  Hence
 \begin{equation}\label{conxlogx}
 \frac{\sum_{j=1}^r t_j}{r} \log \frac{\sum_{j=1}^r t_j}{r}\le
 \frac{1}{r}\sum_{j=1}^r t_j\log t_j, \textrm{ for }
 t_1,\ldots,t_r\in \R_+.
 \end{equation}
 Clearly, the above inequality is equivalent to the inequality
 \begin{equation}\label{conxlogx1}
 (\sum_{j=1}^r t_j)^{\sum_{j=1}^r t_j}\le r^{\sum_{j=1}^r t_j}
 \prod_{j=1}^r t_i^{t_j} \textrm{ for }
 t_1,\ldots,t_r\in \R_+.
 \end{equation}
 Here $0^0=1$.

 \begin{lem}\label{hafexpin}  Let $A=[a_{ij}]\in
 \{0,1\}^{(2n)\times (2n)}\cap \rS_{2n,0}$.  Then for each
 $i\in\an{2n}$
 \begin{equation}\label{hafexpin1}
 (\haf A)^{\haf A}\le r_i^{\haf A} \prod_{j, a_{ij}=1} (\haf
 A(i,j))^{\haf A(i,j)}.
 \end{equation}
 \end{lem}
 \begin{proof}  Let $t_j=\haf A(i,j)$ for $a_{ij}=1$.
 Use (\ref{hafrowex1}) and (\ref{conxlogx1}) to deduce (\ref{hafexpin1}).
 \end{proof}

 To prove our main result we need the following two lemmas.
 \begin{lem}\label{rfacdecseq}  The sequence
 $(k!)^{\frac{1}{k}}, k=1,\ldots,$ is an increasing sequence.
 \end{lem}

 \begin{proof}  Clearly, the inequality $ (k!)^{\frac{1}{k}}<
 ((k+1)!)^{\frac{1}{k+1}}$ is equivalent to the inequality
 $(k!)^{k+1}< ((k+1)!)^{k}$, which is in turn equivalent to
 $k!< (k+1)^k$, which is obvious.
 \end{proof}

 \begin{lem}\label{specseq}  For an integer $r\ge 3$ the
 following inequality holds.
 \begin{equation}\label{specseq1}
 (r!)^{\frac{1}{r}}((r-2)!)^{\frac{1}{r-2}}<((r-1)!)^{\frac{2}{r-1}}.
 \end{equation}
 \end{lem}
 \begin{proof} Raise the both sides of (\ref{specseq1}) to the
 power $r(r-1)(r-2)$ to deduce that (\ref{specseq1}) is equivalent to the
 inequality
 $$(r!)^{(r-1)(r-2)}((r-2)!)^{r(r-1)}<((r-1)!)^{2r(r-2)}.$$
 Use the identities
 \begin{eqnarray*}
 r!=r(r-1)!,\quad (r-1)!=(r-1)(r-2)!,\\
 2r(r-2)=(r-1)(r-2)+r(r-1)-2,\quad r(r-1)-2=(r+1)(r-2)
 \end{eqnarray*}
 to deduce
 that the above inequality is equivalent to
 $$r^{(r-1)(r-2)} ((r-2)!)^2< (r-1)^{(r+1)(r-2)}.$$
 Take the logarithm of the above inequality, divide it by $(r-2)$
 deduce that (\ref{specseq1}) is equivalent to the inequality
 $$(r-1)\log r + \frac{2}{r-2}\log (r-2)!-(r+1)\log (r-1)<0.$$
 This inequality is equivalent to
 \begin{equation}\label{specseq2}
 s_r:=(r-1)\log \frac{r}{r-1} + 2\big(\frac{1}{r-2}\log (r-2)!-\log
 (r-1)\big)<0 \textrm{ for } r\ge 3.
 \end{equation}
 Clearly
 $$(r-1)\log \frac{r}{r-1}=(r-1)\log
 (1+\frac{1}{r-1})<(r-1)\frac{1}{r-1}=1.$$
 Hence (\ref{specseq2}) holds if
 \begin{equation}\label{specseq3}
 \frac{1}{r-2}\log (r-2)!-\log (r-1) <-\frac{1}{2}.
 \end{equation}
 Recall the Stirling's formula \cite[pp. 52]{Fel}
 \begin{equation}\label{stirf}
 \log k!=\frac{1}{2}\log (2\pi k) + k\log k -k
 +\frac{\theta_k}{12k} \textrm{ for some } \theta_k\in (0,1).
 \end{equation}
 Hence
 $$\frac{\log (r-2)!}{r-2}< \frac{\log 2\pi (r-2)}{2(r-2)} +\log(r-2) -1
 +\frac{1}{12(r-2)^2}.$$
 Thus
 $$\frac{1}{r-2}\log (r-2)!-\log (r-1)<\frac{\log 2\pi (r-2)}{2(r-2)}
 +\log\frac{r-2}{r-1} +\frac{1}{12(r-2)^2} -1.$$
 Since $e^x$ is convex, it follows that $1+x\le e^x$.  Hence
 $$\frac{1}{r-2}\log (r-2)!-\log (r-1)<\frac{\log 2\pi
 (r-2)}{2(r-2)}-\frac{1}{r-1}+\frac{1}{12(r-2)^2}-1.$$
 Note that $-\frac{1}{r-1}+\frac{1}{12(r-2)^2}<0$ for $r\ge 3$.
 Therefore
 \begin{equation}\label{specseq4}
 \frac{1}{r-2}\log (r-2)!-\log (r-1)<\frac{\log 2\pi
 (r-2)}{2(r-2)}-1.
 \end{equation}
 Observe next that that the function $\frac{\log 2\pi x}{2x}$ is
 decreasing for $x> \frac{e}{2\pi}$.  Hence the right-hand side
 of (\ref{specseq4}) is a decreasing sequence for
 $r=3,\ldots,$.  Since $\frac{\log 2\pi\cdot 3}{2\cdot
 3}=0.4894$, it follows that the right-hand side of
 (\ref{specseq4}) is less than $-0.51$ for $r\ge 5$.
 Therefore (\ref{specseq2}) holds for $r\ge 5$.
 Since
 $$s_3=\log\frac{9}{16}<0, \quad s_4=\log \frac{128}{243}<0$$
 we deduce the lemma.
 \end{proof}

 The arguments of the Proof of Lemma \ref{specseq}
 yield that $s_r, r=3,\ldots,$ converges to $-1$.  We checked
 the values of this sequence for $r=3,\ldots, 100$, and we found that
 this sequence decreases in this range.  We conjecture that the
 sequence $s_r, r=3,\ldots$ decreases.

 \section{Proof of generalized Bregman-Minc inequality}

 \begin{theo}\label{gbnin}  Let $G=(V,E)$
 be undirected simple graph on an
 even number of vertices.
 Then the inequality (\ref{uppmatbd})
 holds.
 \end{theo}

 \begin{proof}
 We prove (\ref{bmhafin}).
 We use the induction on $n$.  For $n=1$ (\ref{bmhafin}) is
 trivial.  Assume that theorem holds for $n=m-1$.  Let $n=m$.
 It is enough to assume that $\haf A >0$.  In
 particular each $r_i\ge 1$.  If $r_i=1$ for some $i$, then
 by expanding $\haf A$ by the row $i$,
 using the induction hypothesis and Lemma \ref{rfacdecseq},
 we deduce easily the theorem in this case.  Hence we assume
 that $r_i\ge 2$ for each $i\in \an{2n}$.
 Let $G=G(A)=(\an{2n},E)$ be the graph induced by $A$.
 Then $\haf A>0$ is the number of perfect matchings in $G$.
 Denote by $\cM:=\cM(G)\subset\cM(K_{2n})$ the set of all
 perfect matchings in $G$.  Then $\#\cM=\haf A$.  We now follow the
 arguments in the proof of the Bregman-Minc theorem given in
 \cite{Sch78} with the corresponding modifications.

 \begin{eqnarray*}
 (\haf A)^{2n\haf A}{\buildrel (1)\over =}
 \prod_{i=1}^{2n} (\haf A)^{\haf A} {\buildrel (2)\over \le}
 \prod_{i=1}^{2n} \big(r_i^{\haf
 A}\prod_{j, a_{ij}=1} (\haf A(i,j))^{\haf A(i,j)}\big)\\
 {\buildrel (3)\over
 =}\prod_{\alpha\in\cM}\big(\big(\prod_{i=1}^{2n}r_i\big)\big(\prod_{i=1}^{2n}
 \haf A(i,\alpha(i)\big)\big)\\
 {\buildrel (4)\over
 \le}\prod_{\alpha\in\cM}\big(\big( \prod_{i=1}^{2n} r_i\big)\prod_{i=1}^{2n}
 \big( \prod_{j\in\an{2n}\backslash\{i,\alpha(i)\}, a_{ij}=a_{\alpha(i)j}=0}
 (r_j!)^{\frac{1}{2r_j}}\big)\\
 \big(\prod_{j\in\an{2n}\backslash\{i,\alpha(i)\}, a_{ij}+a_{\alpha(i)j}=1}
 ((r_j-1)!)^{\frac{1}{2(r_j-1)}}\big)
 \big(\prod_{j\in\an{2n}\backslash\{i,\alpha(i)\}, a_{ij}+a_{\alpha(i)j}=2}
 ((r_j-2)!)^{\frac{1}{2(r_j-2)}}\big)\big)\\
 {\buildrel (5)\over =}
 \prod_{\alpha\in\cM}\big(\big( \prod_{i=1}^{2n} r_i\big)\prod_{j=1}^{2n}
 \big( \prod_{i\in\an{2n}\backslash\{j,\alpha(j)\},
 a_{ij}=a_{\alpha(i)j=0}}(r_j!)^{\frac{1}{2r_j}}\big)\\
 \big( \prod_{i\in\an{2n}\backslash\{j,\alpha(j)\},
 a_{ij}+a_{\alpha(i)j=1}}((r_j-1)!)^{\frac{1}{2(r_j-1)}}\big)
 \big( \prod_{i\in\an{2n}\backslash\{j,\alpha(j)\},
 a_{ij}+a_{\alpha(i)j=2}}((r_j-2)!)^{\frac{1}{2(r_j-2)}}\big)\big)\\
 {\buildrel (6)\over \le}
 \prod_{\alpha\in\cM}\big(\big( \prod_{i=1}^{2n} r_i\big)\prod_{j=1}^{2n}
 \big((r_j!)^{\frac{2n-2r_j}{2r_j}}\big)\big(((r_j-1)!)^{\frac{2(r_j-1)}{2(r_j-1)}}\big)\big)\\
 {\buildrel (7)\over =} \prod_{\alpha\in\cM}\big(\big( \prod_{i=1}^{2n}
 (r_i!)^{\frac{2n}{2r_i}}\big){\buildrel (8)\over =}\big( \prod_{i=1}^{2n}
 (r_i!)^{\frac{1}{2r_i}}\big)^{2n\haf A}.
 \end{eqnarray*}
 We now explain each step of the proof.
 \begin{enumerate}
 \item Trivial.
 \item Use (\ref{hafexpin1}).
 \item The number of factors of $r_i$ is equal to $\haf A$ on
 both sides, while the number of factors $\haf A(i,j)$
 equals to the number of $\alpha\in\cM$ such that
 $\alpha(i)=j$.
 \item Apply the induction hypothesis to each $\haf A(i,\alpha(i))$.
 Note that since the edge $(i,\alpha(i))$ appears in the
 perfect matching $\alpha\in\cM$, it follows that $\haf
 A(i,\alpha(i)) \ge 1$.  Hence if
 $j\in\an{2n}\backslash\{i,\alpha(i)\}$ and $r_j=2$ we must
 have that $a_{ij}+a_{\alpha(i)j}\le 1$.
 \item Change the order of multiplication.
 \item Fix $\alpha\in\cM$ and $j\in\an{2n}$.  Then $j$ is
 matched with $\alpha(j)$.  Consider all other $n-1$ edges
 $(i,\alpha(i))$ in $\alpha$.  $j$ is connected to $r_j-1$
 vertices in $\an{2n}\backslash\{j,\alpha(j)\}$.
 Assume there are $s$ triangles formed by $j$ and the
 $s$ edges out of $n-1$ edges in $\alpha\backslash
 (j,\alpha(j))$.  Then $j$ is connected to $t=r_j-1-2s$
 edges vertices $i\in\an{2n}\backslash\{j,\alpha(j)\}$ such that $j$
 is not connected to $\alpha(i)$.   Hence there are $2n-2-(2t+2s)$
 vertices $k\in\an{2n}\backslash\{j,\alpha(j)\}$such that $j$ is
 not connected to $k$ and $\alpha(k)$.  Therefore, for this
 $\alpha$ and $j$ we have the following terms in (5):
 \begin{eqnarray}\nonumber
 \big( \prod_{i\in\an{2n}\backslash\{j,\alpha(j)\},
 a_{ij}=a_{\alpha(i)j=0}}(r_j!)^{\frac{1}{2r_j}}\big)
 \big( \prod_{i\in\an{2n}\backslash\{j,\alpha(j)\},
 a_{ij}+a_{\alpha(i)j=1}}((r_j-1)!)^{\frac{1}{2(r_j-1)}}\big)\\
 \nonumber
 \big( \prod_{i\in\an{2n}\backslash\{j,\alpha(j)\},
 a_{ij}+a_{\alpha(i)j=2}}((r_j-2)!)^{\frac{1}{2(r_j-2)}}\big)\big)=\\
 \nonumber
 (r_j!)^{\frac{2n-2-(2s+2t)}{2r_j}}((r_j-1)!)^{\frac{2t}{2(r_j-1)}}
 ((r_j-2)!)^{\frac{2s}{2(r_j-2)!}}=\\
 \label{texpres}
 (r_j!)^{\frac{2n-r_j-1}{2r_j}}((r_j-2)!)^{\frac{r_j-1}{2(r_j-2)!}}
 \big((r_j!)^{-\frac{1}{r_j}}((r_j-2)!)^{-\frac{1}{(r_j-2)}}
 ((r_j-1)!)^{\frac{2}{(r_j-1)}} \big)^{\frac{t}{2}}.
 \end{eqnarray}
 In the last step we used the equality $r_j-1=2s+t$.
 Assume first that $r_j>2$. Use Lemma
 \ref{specseq} to deduce that (\ref{texpres}) increases in $t$.
 Hence the maximum value of (\ref{texpres}) is achieved when
 $s=0$ and $t=r_j-1$.  Then (\ref{texpres}) is equal to
 $$(r_j!)^{\frac{2n-2r_j}{2r_j}}((r_j-1)!)^{\frac{2(r_j-1)}{2(r_j-1)}}.$$
 If $r_j=2$ then, as we explained above, $s=0$.  Hence (\ref{texpres})
 is also equal to the above expression.  Hence (6) holds.
 \item Trivial.
 \item Trivial.
 \end{enumerate}
 Thus
 $$(\haf A)^{2n\haf A}\le \big( \prod_{i=1}^{2n}
 (r_i!)^{\frac{1}{2r_i}}\big)^{2n\haf A}.$$
 This establishes (\ref{bmhafin}).
 \end{proof}

\end{document}